\documentclass[preprint,11pt]{amsart}
\usepackage{amscd,amssymb,amsopn,amsmath,amsthm,graphics,amsfonts,enumerate,verbatim,calc}
\usepackage[dvips]{graphicx}
\usepackage[colorlinks=true,linkcolor=red,citecolor=blue]{hyperref}
\newcommand{\RNum}[1]{\uppercase\expandafter{\romannumeral #1\relax}}

\newcommand{\ra}{\rightarrow}
\newcommand{\norm}{\|}

\newcommand{\supp}{{\rm{Supp\; }}}

\newcommand{\cone}{{\rm{Cone}}}

\newcommand{\se}{{\rm{sec}}}
\newcommand{\im}{{\rm{Im}}}
\newcommand{\dom}{{\rm{dom}}}
\newcommand{\two}{{\rm{II}}}
\newcommand{\li}{\langle}
\newcommand{\ri}{\rangle}

\newtheorem{theorem}{Theorem}[section]
\newtheorem{corollary}[theorem]{Corollary}
\newtheorem{lemma}[theorem]{Lemma}

\newtheorem{proposition}[theorem]{Proposition}

\newtheorem{definition}[theorem]{Definition}
\newtheorem{example}[theorem]{Example}

\theoremstyle{plain}

\newcommand{\re}[1]{{\color{black}{#1}}}

\numberwithin{equation}{section}

\begin{document}
	\title[Local minimality of weak geodesics on prox-regular sets]{Local minimality of weak geodesics on prox-regular subsets of Riemannian manifolds}
	\author[J. Ferrera, M. R. Pouryayevali, and H. Radmanesh]{Juan Ferrera, Mohamad R. Pouryayevali, and Hajar Radmanesh}
	\address{IMI, Departamento de An\'{a}lisis Matem\'{a}tico, Facultad Ciencias Matem\'{a}ticas, Universidad Complutense, 28040, Madrid, Spain}
	\email{ferrera@mat.ucm.es}
	\address{Department of Pure Mathematics, Faculty of Mathematics and Statistics, University of Isfahan, Isfahan,
		81746-73441, Iran}
	\email{pourya@math.ui.ac.ir} \email{h.radmanesh@sci.ui.ac.ir}

	\maketitle
	
	\begin{abstract}
		In this paper we prove that every locally minimizing curve with constant speed in a prox-regular subset of a Riemannian manifold is
		a weak geodesic. Moreover, it is shown that under certain assumptions, every weak geodesic is locally minimizing.
		Furthermore a notion of closed weak geodesics on  prox-regular sets is introduced and a characterization of these curves as nonsmooth critical points of the energy functional is presented.
	\end{abstract}

	\section{Introduction}\label{sec1}
	
	Geodesics and minimizing curves are two basic concepts of Riemannian geometry.  If $M$ is a Riemannian manifold  without boundary, then geodesics are critical points of the energy functionals defined on the Hilbert manifold $H^1(I,M)$ of the admissible curves and they are locally minimizing; see \cite{Klingenberg,lee,Schwartz}. Moreover, they satisfy a second order  differential equation on the manifold and hence the set of all geodesics can be parameterized by position and velocity.
	
	On the other hand,  the study of closed geodesics in Riemannian manifold without boundary is one of the most important issues in  Riemannian geometry and it has played an essential role in dynamical systems; see for instance, \cite{Klingenberg1}. Investigations in closed geodesics  leads to the development of sophisticated tools in nonlinear analysis and symplectic geometry. The standard technique to specify the closed geodesics is critical point theory.
	
	When the manifold $M$ has boundary, even if it is smooth,  or in the case of  $p$-convex  (or $\varphi$-convex) subsets of  $\mathbb{R}^n$ different kind of irregularities may be appeared. In the latter case  using a new definition of geodesics in the framework of Sobolev spaces, a characterization of geodesics in terms of local minimality of energy functional on the set of suitable curves is obtained; see \cite{Canino1,Canino2}.
	
	The class of $p$-convex subsets of  $\mathbb{R}^n$  contains submanifolds (possibly with boundary) of class $C^{1,1}_{loc}$, images under $C^{1,1}_{loc}$-diffeomorphisms of convex sets and even some sets which are not topological manifolds. The problem of existence of  closed geodesics is extended to $p$-convex subsets of $\mathbb{R}^n$ where as mentioned above, in some cases these sets have not a natural structure of manifold, see \cite{Canino}.
	
	There are various equivalent definitions of the notion of $p$-convexity which among them one can mention the concept of prox-regularity of sets (see \cite{Colombo1}) which has significant applications in crowed motion and Moreau sweeping process, see \cite{Maury,Tanwani}.
	
	In \cite{Barani}  the class of  $\varphi$-convex subsets of infinite-dimensional Hadamard manifolds was studied and in \cite{Hosseini} the notion of prox-regularity of sets in  Riemannian manifolds was introduced. It was proved in \cite{convex} that these two classes in the setting of  finite dimensional Riemannian manifolds  coincide.

	In \cite{Minimizing}  the subject of minimizing curves on a prox-regular subset $S$ of a Riemannian manifold $M$ with $C^2$ boundary was considered.  Motivated by \cite{Canino2}  the case where the prox-regular set $S\subseteq M$ does not have a $C^2$ boundary was considered and by introducing weak geodesics on $S$ in the intrinsic case, it was  proved that weak geodesics are critical points of the energy functional, see  \cite{weak}.
	
	Recall that a weak geodesic on $S$ is a curve $\gamma\in W^{2,2}(I,M)$ such that  $\im \gamma\subseteq S$ and $D_t\dot{\gamma}(t)\in N^P_S\left(\gamma(t)\right)$ for almost all $t\in I$, where $N^P_S\left(\gamma(t)\right)$ denotes the proximal normal cone to $S$ at $\gamma(t)$.
	The above definition depends on the chosen normal cone and hence it does not seem to be satisfactory.
	The purpose of this paper is to define an intrinsic characterization of geodesics in terms of local minimality
	for the energy functional. More precisely, we show that every locally minimizing curve with constant speed in the prox-regular set
	$S$ is a weak geodesic and conversely, a weak geodesic $\gamma$ on $S$ is locally minimizing provided that  the sectional curvatures
	of $M$ are nonpositive at least in a neighborhood of $\im(\gamma)$. In order to remove the assumption of
	non-positiveness of sectional curvatures, we first consider the problem on a complete Riemannian submanifold $M$ of $\mathbb{R}^n$ and we
	demonstrate that every prox-regular subset $S$ of $M$ is also a prox-regular subset of $\mathbb{R}^n$. Then we show that every weak geodesic on $S$
	as a subset of $M$ is also a weak geodesic on $S$ as a subset of $\mathbb{R}^n$ and using an isometric embedding,
	we get the same results for complete Riemannian manifolds. This implies that on prox-regular subsets of complete
	Riemannian manifolds, every weak geodesic is locally minimizing. Moreover, we define closed weak geodesics in prox-regular subsets
	of Riemannian manifolds and a characterization of these geodesics is given.
	
	The paper is organized as follows. In Section \ref{sec2}, we recall some relevant definitions and notations from nonsmooth analysis and Riemannian geometry.
	Section \ref{sec3} is concerned with the local minimality property of weak geodesics and we discuss the relation between weak
	geodesics and locally minimizing curves in $S$.
	In Section \ref{sec4}, we obtain a description of prox-regular sets and weak geodesics in complete Riemannian embedded submanifolds of Euclidean spaces.
	Section \ref{sec5} is devoted to  closed weak geodesics on $S$ and we give two characterizations of these notions in terms of nonsmooth critical points of the energy functional and local minimality.

	\section{Preliminaries and notations}\label{sec2}
	Let us first recall some standard notations and known results of Riemannian manifolds and nonsmooth analysis; see, e.g.,
	\cite{Azagra1,Clark book a,Do Carmo,sakai}. Throughout this paper, $(M,g)$ is an $n$-dimensional Riemannian manifold endowed with a Riemannian metric $g_x=\li .,.\ri_x$ on each tangent space $T_xM$ and the Riemannian connection $\nabla$. Moreover, $I$ is used to denote the closed interval $[0,1]$ and  for every $x, y \in M$,  $d(x,y)$ denotes the
	Riemannian distance from $x$ to $y$.
	
	The map $\exp_x : U_x \to M$ is the exponential map at $x$, where $U_x$ is an open subset of the tangent space $T_xM$
	containing $0_{x}$.
	Furthermore for $x\in M$, let $r(x)$ be the convexity radius at $x$, then the function $x\mapsto r(x)$ from
	$M$ to $\mathbb{R}^{+} \cup \{+\infty \}$ is continuous; see \cite{sakai}.
	For a curve $\gamma : I \rightarrow M$ and $t_0,t\in I$, the notation  $L^{\gamma}_{t_0,t}$ signifies the parallel
	transport along $\gamma$ from $\gamma(t_0)$ to $\gamma(t)$. When $\gamma$ is the unique minimizing geodesic
	joining $y$ to $x$, we use the notation $L_{y,x}$.
	
	
	We use the notations  $C^0(I,M)$ and $C^{\infty}(I,M)$ to denote the space of continuous curves, endowed with the metric $d_{\infty}\left(\gamma,\eta\right)=\sup_{t\in I}d\left(\gamma(t),\eta(t)\right)$ and  the set of smooth curves, respectively.
	Let us now recall the Sobolev space of  manifold valued curves.  The following proposition  is a slight modification of \cite[Lemma B.5]{Wehrheim}.
	\begin{proposition}\label{prop1}
		Let $k\in \mathbb{N}$ and $1\leq p<\infty$ be such that $kp>1$. Then for $\gamma\in C^0(I,M)$ the following statements are equivalent:
		\begin{itemize}
			\item[(i)] $\phi_{\alpha}\circ\gamma\in W^{k,p}\left(\gamma^{-1}\left(U_{\alpha}\right), \mathbb{R}^{n}\right)$ for every chart  $\left(U_{\alpha}, \phi_{\alpha}\right)$ on $M$;
			\item[(ii)] $\Phi\circ \gamma\in W^{k,p}\left(I,\mathbb{R}^{2n+1}\right)$;
			\item[(iii)] $\gamma=\exp_c(V)$ for some $c\in C^{\infty}(I,M)$ and $V\in W^{k,p}\left(I, \gamma^{-1}TM\right)$,
		\end{itemize}
		where $\gamma^{-1}TM$ denotes the pullback bundle of $TM$ by $\gamma$  and $\Phi: M\ra \mathbb{R}^{2n+1}$ is an embedding.
	\end{proposition}
		The set of continuous curves that satisfy these equivalent statements is called the Sobolev space  of curves on $M$ and denoted by $W^{k,p}\left(I,M\right)$. For $k=0$ and $p=2$,  the space $L^2(I,M)$ is defined as
	\[
	L^2(I,M):=\left\{\gamma:I\ra M : \Phi\circ \gamma\in L^2\left(I,\mathbb{R}^{2n+1}\right)\right\},
	\]
	with the topology given by
	\[
	\gamma_i\ra \gamma \ \text{in}\ L^2(I,M)\quad\Leftrightarrow\quad \Phi\circ\gamma_i\ra \Phi\circ\gamma \ \text{in}\ L^2\left(I,\mathbb{R}^{2n+1}\right).
	\]
	Moreover, the Hilbert manifold  $W^{1,2}\left(I,M\right)$ is denoted by $H^1(I,M)$. For more details regarding Sobolev spaces, see \cite{Adams,Convent1,Convent2}.

	For  a nonempty closed subset $S$ of  $M$, the distance function to $S$ is $d_S(z)= \inf_{x\in S} \, d(x,z)$    and
	\[
	P_S(z)= \left\{ x\in S : d_S(z)=d(x,z)\right \}\quad \forall z\in M,
	\]
is the metric projection to $S$. Moreover, $N^{P}_{S}(x)$ denotes the proximal normal cone to $S$ at $x\in S$  and $\xi\in N^{P}_{S}(x)$ if and only if there exist  $\sigma>0$ and a convex neighborhood $U$ of $x$ such that
	\[
	\li \xi, \exp_x^{-1}y \ri \leq\sigma\;d^2(x,y),
	\]
	for every $y\in U\cap S$.
	
	For a lower semicontinuous function $f:M\ra (-\infty,+\infty]$   and $x\in \dom(f)$,  the set
	\[
	D^-f(x):=\left\{dg(x): g\in C^1(M,\mathbb{R}), f-g \ \text{attains a local minimum at}\ x\right\},
	\]
	is called the viscosity (or Fr\'{e}chet) subdifferential of $f$ at $x$.
	
	A closed subset $S$ of $M$ is called prox-regular at $\bar{x}\in S$ if there exist $\varepsilon>0$ and $\sigma>0$ such that  $B(\bar {x}, \varepsilon)$ is convex and $\li v, \exp_x^{-1}y \ri \leq\sigma d^2(x,y)$ 	
for every  $x,y\in S\cap B(\bar{x}, \varepsilon)$ and $v\in N^P_S(x)$ with $\norm v\norm <\epsilon$.
	Moreover, $S$ is  prox-regular if it is prox-regular at each point of $S$; for more details, see \cite{Hosseini}. 	
	 Every prox-regular subset $S$ of $M$ is  $\varphi$-convex for a continuous function $\varphi : S\rightarrow [0,\infty)$  (\cite{convex}).
	Recall that a nonempty closed set $S\subset M$ is called $\varphi$-convex if for every $x\in S$ and $v\in N^P_S(x)$
	\[
	\li v,\exp_x^{-1}y \ri \leq \varphi(x)\norm v\norm d^2(x,y)  \quad \ \forall  y\in U\cap S,
	\]
	 where $U$ is a convex neighborhood of $x$.
	In \cite{Minimizing} it is shown that for a closed $\varphi$-convex set $S$, the metric projection $P_S$ is locally Lipschitz on an open set containing  $S$ and it is directionally differentiable at each point $x\in S$.


	\section{Local minimality of weak geodesics}\label{sec3}
	In this section, by defining locally minimizing curves in a prox-regular set $S$ in terms of length functional, we get another equivalent definition in terms of energy functional. Then motivated by \cite{Canino1}, we show that a locally minimizing curve with constant speed in $S$ is a weak geodesic. Finally, we prove that the converse holds provided that the sectional curvatures of $M$ are nonpositive in a neighborhood of  the image of a weak geodesic.

	Let $S$ be a  closed and connected  prox-regular subset of $M$ and $U$ be an open neighborhood of $S$ on which $P_S$ is single-valued and locally Lipschitz. Moreover, let $\varphi : S\rightarrow [0,\infty)$ be a continuous function such that $S$ is  $\varphi$-convex. For any $x,y\in S$, the set of admissible curves $\mathcal{A}_{x,y}$ is considered to be
	\[
	\mathcal{A}=\mathcal{A}_{x,y}:=\left\{\eta\in H^1(I,M) : \eta(t)\in S, \ \forall t\in I, \eta(0)=x, \eta(1)=y\right\}.
	\]
	We now recall the definition of a weak geodesic on $S$ that has been introduced in \cite{weak}.
	\begin{definition}
		A continuous curve $\gamma : I\ra M$  is called a weak geodesic on $S$ joining $x,y\in S$ if $\gamma(0)=x$, $\gamma(1)=y$ and
		\begin{itemize}
			\item [(a)] $\gamma(t)\in S \quad \forall t\in I$;
			\item [(b)] $\gamma\in W^{2,2}(I,M)$;
			\item [(c)] $D_t\dot{\gamma}(t)\in N^P_S\left(\gamma(t)\right) \quad a.\ e.\ t\in I$.
		\end{itemize}
	\end{definition}
	
	In \cite{weak} we characterized weak geodesics on $S$ as viscosity critical  points  of the energy functional defined on  $L^2(I,M)$ by
	\[
	f(\eta)=f_{x,y}(\eta):=\left\{
	\begin{array}{lr}
		\frac{1}{2} \int_0^1 \norm \dot{\eta}(t)\norm ^2 dt & \eta\in \mathcal{A} \\
		& \\
		+\infty & \eta\in L^2(I,M)\setminus\mathcal{A}.
	\end{array}\right.
	\]
	Indeed, we proved that the  set $\mathcal{A}$ is closed in $L^2(I,M)$ and the energy functional $f$ is lower semicontinuous and finally, 
	we proved that $\gamma$ is a weak geodesic on $S$ if and only if $0\in D^-f(\gamma)$. Moreover, if $\gamma\in \mathcal{A}$ is a weak geodesic on $S$, then $\gamma\in W^{2,\infty}\left(I, M\right)$ and $\gamma$ has constant speed.
	
	Note that for any $\eta\in H^1(I,M)$ we have $\mathcal{L}\left(\eta\right)\leq \sqrt{2f\left(\eta\right)}$, where $\mathcal{L}\left(\eta\right)$ denotes the length of $\eta$ defined by
	\[
	\mathcal{L}\left(\eta\right):=\int_{0}^{1}\norm \dot{\eta}(t)\norm dt,
	\]
	and the equality holds if and only if $\norm \dot{\eta}(t)\norm$ is constant for almost all $t\in I$.
	
	Moreover, we recall that for any absolutely continuous curve $\eta : I\ra M$, the derivative $\dot{\eta}(t)$ exists for almost all $t\in I$ and $\norm \dot{\eta}\norm\in L^1(I)$; see \cite{Burtscher}. Hence   $\mathcal{L}\left(\eta\right)$ is well-defined  and  using \cite[Lemma 3.14]{Burtscher} and \cite[Theorem 3.15]{Burtscher}, the function $t\mapsto \mathcal{L}\left(\eta|_{[0,t]}\right)$ is also  absolutely continuous on $I$. Furthermore according to \cite[Corollary 3.8]{Burtscher}
	\begin{equation}\label{derac}
		\lim_{h\ra 0}\frac{d\left(\eta(t),\eta(t+h)\right)}{|h|}=\|\dot{\eta}(t)\|,  \quad a.e.\ t\in I.
	\end{equation}

	
	The class of  absolutely continuous curves on $M$ is denoted by $AC\left(I,M\right)$ and for any $x,y\in S$ we consider the following subclass  of $AC\left(I,M\right)$ as
	\[
	\mathcal{B}=\mathcal{B}_{x,y}:=\left\{\eta\in AC\left(I,M\right) : \im\left(\eta\right)\subset S, \ \eta(0)=x, \ \eta(1)=y\right\}.
	\]
	We are now ready to define locally minimizing curves on $S$.
	
	\begin{definition}
		A curve $\gamma\in AC(I,M)$ with $\im (\gamma)\subseteq S$ is said to be locally minimizing, if for every $t_0\in I$ there exists $h>0$ such that $\mathcal{L}(\gamma)\leq \mathcal{L}(\eta)$ for every  $\eta\in AC(I,M)$ with  $\im (\eta)\subseteq S$ and $\eta=\gamma$ on $I\setminus \left(t_0-h,t_0+h\right)$.
	\end{definition}
	
	Let us now obtain a characterization of locally minimizing curves on $S$ in terms of the energy functional.
	To this end, we show that a curve $\gamma\in \mathcal{B}$ with constant speed minimizes the length functional $\mathcal{L}$ on $\mathcal{B}$ if and only if $\gamma$  minimizes the energy functional $f$ on $L^2(I,M)$.
	
	\begin{lemma}\label{repa}
		Let $\eta:I\ra M$ be an absolutely continuous curve. Then there exists a nondecreasing  absolutely continuous function $\chi:I\ra I$  with $\chi(0)=0$, $\chi(1)=1$ and a Lipschitz curve $\tilde{\eta}:I\ra M$   with the properties that $\eta=\tilde{\eta} \circ \chi$, $\|\dot{\tilde{\eta}}(t)\|=\mathcal{L}(\eta)$ for almost all $t\in I$ and $\mathcal{L}(\tilde{\eta})=\mathcal{L}(\eta)$.
	\end{lemma}
	\begin{proof}
		Without loss of generality, we assume that $\eta$ is not constant. For every $t\in I$, we define
		\[
		\chi(t):= \frac{1}{\mathcal{L}(\eta)}\mathcal{L}\left(\eta|_{[0,t]}\right),
		\]
		then clearly $\chi$ is nondecreasing and $\chi(0)=0$, $\chi(1)=1$. Moreover, $\chi$ is absolutely continuous and hence is onto.
		
		Since $\eta$ is constant on the set $\chi^{-1}(t)$ for every $t\in I$, the curve $\tilde{\eta}:I\ra M$ defined by $\tilde{\eta}(t):=\eta \circ \chi^{-1}(t)$ is well-defined.
		We show that $\tilde{\eta}$ is Lipschitz. Indeed, for $s_1,s_2\in I$ we have
		\begin{align*}
			d\left(\tilde{\eta}\left(s_1\right),\tilde{\eta}\left(s_2\right)\right) =& d\left(\tilde{\eta}\left(\chi(t_1)\right), \tilde{\eta}\left(\chi(t_2)\right)\right)\\
			= & d\left(\eta(t_1),\eta(t_2)\right)\leq \int_{t_1}^{t_2} \|\dot{\eta}(t)\|dt\\
			= & \mathcal{L}(\eta)\left(\chi(t_2)-\chi(t_1)\right) \\
			\leq & \mathcal{L}(\eta)|s_1-s_2|,
		\end{align*}
		where $t_1,t_2\in I$ are the points such that $s_i=\chi(t_i)$ for $i=1,2$.
		
		Furthermore by \eqref{derac} for almost all $t\in I$ we have
		\begin{align*}
			\|\dot{\tilde{\eta}}(t)\|= & \lim_{h\ra 0}\frac{d\left(\tilde{\eta}(t),\tilde{\eta}(t+h)\right)}{|h|} \\
			\leq & \lim_{h\ra 0}\frac{\mathcal{L}(\eta)|h|}{|h|}=\mathcal{L}(\eta),
		\end{align*}
		and hence $\|\dot{\tilde{\eta}}(t)\|\leq \mathcal{L}(\eta)$ for almost all $t\in I$. On the other hand, by a change of variable we obtain that
		\[
		\int_{0}^{1}\|\dot{\tilde{\eta}}(t)\|dt=\int_{0}^{1}\|\dot{\eta}(t)\|dt=\mathcal{L}(\eta),
		\]
		and we conclude that $\|\dot{\tilde{\eta}}(t)\|=\mathcal{L}(\eta)$, for almost all $t\in I$.
	\end{proof}
	
	\begin{theorem}\label{fl}
		Suppose that $\gamma\in \mathcal{B}$, then the following statements are equivalent:
		\begin{itemize}
			\item[(i)] $\mathcal{L}(\gamma)\leq \mathcal{L}(\eta)$, for all $\eta\in \mathcal{B}$ and $\norm \dot{\gamma}\norm$ is constant a.e. on $I$;
			\item[(ii)] $f(\gamma)\leq f(\eta)$, for all $\eta\in \mathcal{A}$.
			
		\end{itemize}
	\end{theorem}
	\begin{proof}
		Using Lemma \ref{repa} 
		we get the result.
	\end{proof}
	Using Theorem \ref{fl}, we can characterize  locally minimizing curves on $S$ using the energy functional. This characterization is more applicable to study the local minimality for weak geodesics.
	\begin{theorem}\label{el}
		Let $\gamma \in \mathcal{B}$, then the following statements are equivalent:
		\begin{itemize}
			\item[(i)]  $\gamma$ is locally minimizing and $\norm \dot{\gamma}\norm$ is constant a.e. on $I$;
			\item[(ii)] for every $t_0\in I$ there exists $h>0$ such that $f(\gamma)\leq f(\eta)$ for every  $\eta\in \mathcal{A}$ with the property that $\eta=\gamma$ on $I\setminus \left(t_0-h,t_0+h\right)$.
		\end{itemize}
	\end{theorem}
	
	We now show that every locally minimizing curve is a weak geodesic. To this end, we need to prove the following lemmas.
	
	\begin{lemma}\label{dt}
		Let $\gamma:I\ra M$ be a smooth curve and $s\in (0,1)$ be fixed. Then
		\[
		D_t\dot{\gamma}(s)=\lim_{h\ra 0}\frac{\exp_{\gamma(s)}^{-1}\gamma(s+h)+\exp_{\gamma(s)}^{-1}\gamma(s-h)}{h^2}.
		\]
	\end{lemma}
	\begin{proof}
		Let  $h>0$ be small enough such that the image of $\gamma$ on $(s-h,s+h)$ is contained entirely  in a convex ball around $\gamma(s)$. We now define
		\[
		g(t):=\exp_{\gamma(s)}^{-1}\gamma(t), \quad \forall t\in (s-h,s+h).
		\]
		Then we have
		\[
		\hspace{-4cm}  \lim_{h\ra 0}\frac{\exp_{\gamma(s)}^{-1}\gamma(s+h)+ \exp_{\gamma(s)}^{-1}\gamma(s-h)}{h^2}
		\]
		\begin{align*}
			&=\lim_{h\ra 0}\frac{g(s+h)+g(s-h)-2g(s)}{h^2} \\
			& =\frac{d^2}{dt^2}|_{t=s}g(t) \\
			& =\frac{d}{dt}|_{t=s}\left\li d\exp_{\gamma(s)}^{-1}\left(\gamma(t)\right),\dot{\gamma}(t)\right\ri \\
			& =\left\li d^2\exp_{\gamma(s)}^{-1}\left(\gamma(s)\right)\left(\dot{\gamma}(s)\right),\dot{\gamma}(s)\right\ri\\
			& \qquad \qquad\qquad +\left\li d\exp_{\gamma(s)}^{-1}\left(\gamma(s)\right),D_t\dot{\gamma}(s)\right\ri\\
			& = D_t\dot{\gamma}(s),
		\end{align*}
		since for every $p\in M$ and $v\in T_pM$, we have $d^2\exp_p^{-1}(p)(v,v)=0$.
	\end{proof}

	\begin{lemma}\label{dtnorm}
		Let $\gamma:I\ra M$ be a differentiable curve such that $\im(\gamma)\subseteq S$ and $s\in (0,1)$. If $D_t\dot{\gamma}(s)$ exists, then
		\[
		\left \norm \mathrm{P}_{N^P_S\left(\gamma(s)\right)}\left(D_t\dot{\gamma}(s)\right)\right\norm\leq 2\varphi\left(\gamma(s)\right)\norm \dot{\gamma}(s)\norm^2,
		\]
		where $\mathrm{P}_{N^P_S\left(\gamma(s)\right)}$ denotes the projection on the proximal normal cone $N^P_S\left(\gamma(s)\right)$.
	\end{lemma}
	\begin{proof}
		Let $v\in N^P_S\left(\gamma(s)\right)$. Since $\im(\gamma)\subseteq S$, for all $h>0$ small enough we have
		\[
		\left\li v, \exp_{\gamma(s)}^{-1}\gamma(s+h)\right\ri\leq \varphi\left(\gamma(s)\right)\norm v\norm d^2\left(\gamma(s),\gamma(s+h)\right),
		\]
		and
		\[
		\left\li v, \exp_{\gamma(s)}^{-1}\gamma(s-h)\right\ri\leq \varphi\left(\gamma(s)\right)\norm v\norm d^2\left(\gamma(s),\gamma(s-h)\right).
		\]
		Therefore
		\[
		\hspace{-4cm}\left\li v, \frac{\exp_{\gamma(s)}^{-1}\gamma(s+h)+ \exp_{\gamma(s)}^{-1}\gamma(s-h)}{h^2}\right\ri\leq
		\]
		\[
		\varphi\left(\gamma(s)\right)\norm v\norm\left(\left(\frac{d\left(\gamma(s),\gamma(s+h)\right)}{h}\right)^2+
		\left(\frac{d\left(\gamma(s),\gamma(s-h)\right)}{h}\right)^2\right).
		\]
		 Taking the limit when $h\ra 0$ and using Lemma \ref{dt}, we obtain that
		\[
		\left\li v, D_t\dot{\gamma}(s)\right\ri\leq 2\varphi\left(\gamma(s)\right)\norm v\norm\norm \dot{\gamma}(s)\norm^2,
		\]
		for all $v\in N^P_S\left(\gamma(s)\right)$. Hence putting $v:=\mathrm{P}_{N^P_S\left(\gamma(s)\right)}\left(D_t\dot{\gamma}(s)\right)$ and noting that $N^P_S\left(\gamma(s)\right)$ is closed and convex, we get the result.
	\end{proof}
	
	\begin{theorem}\label{localweak}
		Every locally minimizing curve in $S$ is a weak geodesic when it is given a constant speed parametrization.
	\end{theorem}
	\begin{proof}
		Let $\gamma\in H^1(I,M)$ with $\im (\gamma)\subseteq S$ be a locally minimizing curve with constant speed and $t_0\in I$. Then according to Theorem \ref{el}, there exists $h>0$ such that $f(\gamma)\leq f(\eta)$ for every  $\eta\in \mathcal{A}$ with the property that $\eta=\gamma$ on $I\setminus \left(t_0-h,t_0+h\right)$.
		
		For convenience  we assume that $t_0\in (0,1)$. We consider the restricted curve $\tilde{\gamma}:=\gamma|_{[t_0-h,t_0+h]}$ and we put $I_0:=[t_0-h,t_0+h]$, $\gamma(t_0-h):=x_0$ and $\gamma(t_0+h):=y_0$. Then we have $f_{x_0,y_0}\left(\tilde{\gamma}\right)\leq f_{x_0,y_0}(\eta)$ for all $\eta\in H^1(I_0,M)$ with $\im(\eta)\subseteq S$ and $\eta(t_0-h)=x_0$, $\eta(t_0+h)=y_0$. This implies that $0\in D^-f_{x_0,y_0}\left(\tilde{\gamma}\right)$ and hence according to \cite[Theorem 4.3]{weak} the curve $\tilde{\gamma}$ is a weak geodesic on $S$ joining $x_0,y_0$. Indeed,  $\gamma\in W^{2,2}\left(I_0,M\right)$ and
		\[
		D_t\dot{\gamma}(t)\in N^P_S\left(\gamma(t)\right) \quad a.\ e.\ t\in I_0.
		\]
		Moreover, by \cite[Proposition 4.4]{weak}, $\gamma\in W^{2,\infty}\left(I_0,M\right)$ and $\norm \dot{\gamma}\norm$ is constant on $I_0$.
		Since $t_0\in I$ is arbitrary, we have $\gamma\in W^{2,2}_{loc}\left(I,M\right)$ and
		\[
		D_t\dot{\gamma}(t)\in N^P_S\left(\gamma(t)\right) \quad a.\ e.\ t\in I.
		\]
		
		It suffices to show that $\gamma\in W^{2,2}\left(I,M\right)$. According to \cite[Theorem B.2]{Wehrheim} $\gamma$ is $C^1$ on $I$ and hence $\norm \dot{\gamma}\norm$ as a continuous locally constant function on the connected interval $I$, is constant on $I$. Then $\dot{\gamma}\in L^{\infty}\left(I, \gamma^{-1}TM\right)$ and hence using Lemma \ref{dtnorm} we conclude that
		\[
		\norm  D_t\dot{\gamma}(t)\norm \leq 2\bar{\varphi}\norm \dot{\gamma}\norm^2_{L^{\infty}}  \quad a.e. \ t\in I,
		\]
		where $\bar{\varphi}:=\max_{I}\varphi\circ\gamma$ and  the proof is complete.
	\end{proof}
	
	To justify the converse statement of Theorem \ref{localweak}, we first ask if every weak geodesic on $S$ is locally minimizing.
	Note that in \cite[Theorem 2.7]{Canino1} it is proved that every weak geodesic on a $\varphi$-convex subset of a Hilbert space is locally minimizing.
	
	It is worth mentioning that if $g\in W^{1,2}\left([a,b],\mathbb{R}\right)$ and $g(a)=g(b)=0$, then
	\begin{equation}\label{ineqq}
		\sup_{t\in [a,b]}|g(t)|\leq \left(\frac{b-a}{2}\right)^{1/2}\norm g'\norm_{L^2}.
	\end{equation}

	
	\begin{lemma}\label{dinfty}
		Let $\gamma,\eta\in H^1(I,M)$ and there exist $t_0\in I$ and $h>0$ such that $\eta=\gamma$ on $I\setminus (t_0-h,t_0+h)$. Then
		\[
		d_{\infty}^2\left(\gamma,\eta\right)\leq h\int_{0}^{1}\norm L_{\eta,\gamma}\dot{\eta}-\dot{\gamma}\norm^2dt,
		\]
		provided that $d_{\infty}\left(\gamma, \eta\right)< \bar{r}$, where $\bar{r}=\min_{t\in I}r\left(\gamma(t)\right)$.
	\end{lemma}
	\begin{proof}
		We define $g(t):=d\left(\gamma(t),\eta(t)\right)$ for all $t\in I$. Hence $g\equiv0$ on $I\setminus (t_0-h,t_0+h)$ and $g\in W^{1,2}\left(I,\mathbb{R}\right)$. On the other hand, for almost all $t\in (t_0-h,t_0+h)$ we have
		\begin{align*}
			g'(t)= & \frac{\partial d}{\partial x}\left(\gamma(t),\eta(t)\right) \left(\dot{\gamma}(t)\right)+\frac{\partial d}{\partial y}\left(\gamma(t),\eta(t)\right) \left(\dot{\eta}(t)\right) \\
			= & \left\li -\frac{\exp_{\gamma(t)}^{-1}\eta(t)}{\norm\exp_{\gamma(t)}^{-1}\eta(t) \norm},\dot{\gamma}(t)\right\ri+ \left\li -\frac{\exp_{\eta(t)}^{-1}\gamma(t)}{\norm\exp_{\eta(t)}^{-1}\gamma(t) \norm},\dot{\eta}(t)\right\ri\\
			= & \left\li \frac{\exp_{\gamma(t)}^{-1}\eta(t)}{\norm\exp_{\gamma(t)}^{-1}\eta(t) \norm},L_{\eta(t),\gamma(t)}\dot{\eta}(t)-\dot{\gamma}(t)\right\ri.
		\end{align*}
		Hence
		\[
		|g'(t)|\leq \norm L_{\eta(t),\gamma(t)}\dot{\eta}(t)-\dot{\gamma}(t) \norm  \quad a.e.\  t\in I.
		\]
		Therefore by \eqref{ineqq}, we conclude that
		\[
		\left(\sup_{t\in [t_0-h,t_0+h]} d\left(\gamma(t),\eta(t)\right)\right)^2\leq h \int_{0}^{1}\norm L_{\eta,\gamma}\dot{\eta}-\dot{\gamma}\norm^2dt.
		\]
	\end{proof}
	
	\begin{theorem}\label{weaklocal}
		Let $\gamma\in \mathcal{A}$ be a weak geodesic on $S$. If in a neighborhood of $\im(\gamma)$ the sectional curvatures of $M$ are nonpositive, then $\gamma$ is locally minimizing.
	\end{theorem}
	\begin{proof}
		Let $W$ be a tubular  neighborhood of $\im(\gamma)$ with compact closure and with radius smaller than $\bar{r}$, where $\bar{r}=\min_{t\in I}r\left(\gamma(t)\right)$. We first suppose that $\delta\leq\se\leq\Delta$ on $W$, $\Delta\geq 0$.  Note that for any $\eta\in \mathcal{A}$ with the property that $\im(\eta)\subset W$, we have
		\[
		\norm L_{\eta,\gamma}\dot{\eta}-\dot{\gamma}\norm^2=\norm \dot{\eta}\norm^2+\norm \dot{\gamma}\norm^2-2\left\li \dot{\gamma},L_{\eta,\gamma}\dot{\eta} \right\ri,
		\]
		and so
		\[
		\frac{1}{2}\norm \dot{\eta}\norm^2-\frac{c\left(\gamma,\eta\right)-1}{2}\norm \dot{\gamma}\norm^2=\frac{1}{2}\norm L_{\eta,\gamma}\dot{\eta}-\dot{\gamma}\norm^2+\left\li \dot{\gamma},L_{\eta,\gamma}\dot{\eta}-\frac{c\left(\gamma,\eta\right)}{2}\dot{\gamma} \right\ri,
		\]
		where $c\left(\gamma,\eta\right)(t)=2\sqrt{\Delta}\:d\left(\gamma(t),\eta(t)\right)
		\cot\left(\sqrt{\Delta}\:d\left(\gamma(t),\eta(t)\right)\right)$ for all $t\in I$.
		Moreover, using \cite[Lemma 3.9]{weak} we have
		\[
		\left\li \dot{\gamma}, L_{\eta,\gamma}\dot{\eta}-\frac{1}{2}c\left(\gamma,\eta\right)\dot{\gamma}\right\ri\geq \frac{d}{dt}\left\li \dot{\gamma}, \exp_{\gamma}^{-1}\eta\right\ri-\left\li D_t\dot{\gamma}, \exp_{\gamma}^{-1}\eta\right\ri
		\]
		\[
		- \frac{1}{2}|R|_{\infty} d^2\left(\gamma, \eta\right)\norm \dot{\gamma}\norm \norm \dot{\eta}\norm,
		\]
		where $R$ denotes the curvature tensor on $M$ and thus
		\begin{equation}\label{ineq1}
			\frac{1}{2}\norm \dot{\eta}\norm^2-\frac{c\left(\gamma,\eta\right)-1}{2}\norm \dot{\gamma}\norm^2\geq \frac{1}{2}\norm L_{\eta,\gamma}\dot{\eta}-\dot{\gamma}\norm^2+\frac{d}{dt}\left\li \dot{\gamma}, \exp_{\gamma}^{-1}\eta\right\ri
		\end{equation}
		\[
		-  \left\li D_t\dot{\gamma}, \exp_{\gamma}^{-1}\eta\right\ri-\frac{1}{2}|R|_{\infty} d^2\left(\gamma, \eta\right)\norm \dot{\gamma}\norm \norm \dot{\eta}\norm.
		\]
		Hence by integrating from both side of \eqref{ineq1} and noting that $D_t\dot{\gamma}\in N^P_S\left(\gamma\right),\ a.\ e.$, we have
		\begin{equation}\label{ineq2}
			\frac{1}{2}\int_{0}^{1}\norm \dot{\eta}(t)\norm^2dt -  \frac{1}{2}\int_{0}^{1}\left(c\left(\gamma,\eta\right)-1\right)\norm \dot{\gamma}(t)\norm^2dt\geq
		\end{equation}
		\[
		\frac{1}{2}\int_{0}^{1}\norm L_{\eta,\gamma}\dot{\eta}-\dot{\gamma}\norm^2dt-\int_{0}^{1}\left\li D_t\dot{\gamma},\exp_{\gamma}^{-1}\eta\right\ri dt-\frac{1}{2}|R|_{\infty}\int_{0}^{1} d^2\left(\gamma, \eta\right)\norm \dot{\gamma}\norm \norm \dot{\eta}\norm dt
		\]
		\[
		\geq \frac{1}{2}\int_{0}^{1}\norm L_{\eta,\gamma}\dot{\eta}-\dot{\gamma}\norm^2dt-\int_{0}^{1}\varphi\left(\gamma\right)\norm
		D_t\dot{\gamma}\norm d^2\left(\gamma,\eta\right) dt
		\]
		\[
		-\frac{1}{2}|R|_{\infty}\int_{0}^{1} \norm \dot{\gamma}\norm \norm \dot{\eta}\norm d^2\left(\gamma,\eta\right) dt
		\]
		\[
		\geq \frac{1}{2}\int_{0}^{1}\norm L_{\eta,\gamma}\dot{\eta}-\dot{\gamma}\norm^2dt-\left(\bar{\varphi}\int_{0}^{1}\norm D_t\dot{\gamma}\norm dt+\frac{1}{2}|R|_{\infty} \norm \dot{\gamma}\norm_{L^2} \norm \dot{\eta}\norm_{L^2} \right)d^2_{\infty}\left(\gamma,\eta\right).
		\]
		
		Let $C>0$ be such that $C>\norm \dot{\gamma}\norm_{L^2}$ and $t_0\in (0,1)$. Put
		\[
		\theta:=\bar{\varphi}\int_{0}^{1}\norm D_t\dot{\gamma}(t)\norm dt+\frac{1}{2}C|R|_{\infty}\norm \dot{\gamma}\norm_{L^2},
		\]
		and choose $h>0$ small enough such that $h<\frac{1}{2\theta}$. We claim that for every  $\eta\in \mathcal{A}\setminus\{\gamma\}$ with the properties that $\eta=\gamma$ on $I\setminus \left(t_0-h,t_0+h\right)$, $\norm \dot{\eta}\norm_{L^2}<C$ and $\im(\eta)\subset W$, we have
		\begin{equation}\label{ineq3}
			\frac{1}{2}\int_{0}^{1}\norm \dot{\eta}(t)\norm^2dt > \frac{1}{2}\int_{0}^{1}\left(c\left(\gamma,\eta\right)-1\right)\norm \dot{\gamma}(t)\norm^2dt.
		\end{equation}
		Indeed, from \eqref{ineq2} and using Lemma \ref{dinfty} we obtain that
		\[
		\hspace{-3cm}\frac{1}{2}\int_{0}^{1}\norm \dot{\eta}(t)\norm^2dt- \frac{1}{2}\int_{0}^{1}\left(c\left(\gamma,\eta\right)-1\right)\norm \dot{\gamma}(t)\norm^2dt\geq
		\]
		\[
		\hspace{4cm}\left(\frac{1}{2}-h\theta\right)\int_{0}^{1}\norm L_{\eta,\gamma}\dot{\eta}-\dot{\gamma}\norm^2dt.
		\]
		Moreover, by Lemma \ref{dinfty} we have $\int_{0}^{1}\norm L_{\eta,\gamma}\dot{\eta}-\dot{\gamma}\norm^2dt=0$ if and only if $\eta=\gamma$. It follows that the inequality \eqref{ineq3} holds.
		Note that our claim also holds for $t_0=0$ or $t_0=1$. Indeed, by our choice of $h$, we have $\frac{1-h\theta}{2}>0$.
		
		We now consider the special case when $\Delta=0$ on $W$ and thus $c\left(\gamma,\eta\right)\equiv 2$. In this case, from \eqref{ineq3} we conclude that for every $t_0\in I$ there exists $h>0$ such that $f(\gamma)< f(\eta)$ for every  $\eta\in \mathcal{A}\setminus \{\gamma\}$ with the property that $\eta=\gamma$ on $I\setminus \left(t_0-h,t_0+h\right)$. We note that if $\im(\eta)\nsubseteq W$, then $\mathcal{L}(\gamma)<\mathcal{L}(\eta)$ and it completes the proof.
	\end{proof}

	\begin{example}
		Let $H^2$ denote the hyperbolic plane, that is, the upper half-plane in $\mathbb{R}^2$
		with the metric $g_H=\left(dx^2+dy^2\right)/y^2$. It is well known that  $\left(H^2,g_H\right)$ is a Hadamard manifold and
		for $z_1,z_2\in H^2$ we have
		\[
		d\left(z_1,z_2\right)= 2\ln \frac{\sqrt{\left(x_2-x_1\right)^2+\left(y_2-y_1\right)^2}+\sqrt{\left(x_2-x_1\right)^2+
				\left(y_2+y_1\right)^2}}{2\sqrt{y_1y_2}}.
		\]
		where $z_1=\left(x_1,y_1\right), z_2=\left(x_2,y_2\right)$.
		
		In \cite{Minimizing} we showed that the set $S=\left\{(x,y):1\leq y\leq 2\right\}$ is a prox-regular subset of $H^2$ and
		for every $z=(s,2)\in \partial S$ we have
		\[
		N^P_S(z)=\{\lambda\:\partial/\partial y:\lambda\geq 0\}.
		\]
		For the unit speed curve $\gamma$ in $S$ defined by
		\[
		\gamma(t):=\left(2t,2\right) \qquad \forall t\in \mathbb{R},
		\]
		we have
		$D_t\dot{\gamma}\left(t\right)=2\frac{\partial}{\partial y}$ for all $t\in \mathbb{R}$. It follows that $\gamma$ is a weak geodesic
		on $S$ and hence it is locally minimizing in $S$.
	\end{example}

	Here the question arises whether Theorem \ref{weaklocal} holds in general case without considering the assumption that
	the sectional curvatures of $M$ are nonpositive in a neighborhood of $\im(\gamma)$. In the next section, we will answer
	this question in the case of complete Riemannian manifolds.
	
	\section{Prox-regular subsets of complete Riemannian manifolds and their weak geodesics}\label{sec4}
	
	In this section, we first suppose that $M \subset \mathbb{R}^n$ is a complete Riemannian submanifold of $\mathbb{R}^n$ (at least $C^2$) with the induced metric and we demonstrate  that every prox-regular subset $S$ of $M$ is also a prox-regular subset of $\mathbb{R}^n$.
	Moreover, the relationship between the corresponding proximal normal cones implies that every weak geodesic on $S$ as a subset of $M$
	is also a weak geodesic on $S$ as a subset of  $\mathbb{R}^n$.  We now suppose that $S \subset M$ be a nonempty closed set and $x\in S$. Let us denote the normal space to $M$ at $x$ by $N_x$.
	Moreover, we use the notations $\mathbf{N}^P_S(x)$ and $\mathbf{B}(x,\eta)$ for the proximal normal  cone to $S$ at $x$ by
	considering $S$ as a closed subset of $\mathbb{R}^n$ and the ball in $\mathbb{R}^n$ centered at $x$, respectively.
	
	\begin{lemma}\label{lem4.1}
		Let $M$ be a $C^1$ $m$-dimensional submanifold of $\mathbb{R}^n$ with the induced metric and $x_0\in M$. Then there exist a positive constant $K$ and an open neighborhood $U\subset\mathbb{R}^n$ such that
		\[
		d(x,y)\leq K \norm x-y\norm \qquad \forall x,y\in U\cap M,
		\]
		where $d$ denotes the induced distance function on $M$.
	\end{lemma}
	\begin{proof}
		We may assume, by reordering the coordinates if necessary, that there exists $U\subset\mathbb{R}^n$ such that
		\[
		M\cap U = \left\{\left(\tilde{x}, f(\tilde{x})\right) : \tilde{x}\in\tilde{U}\right\},
		\]
		where $\tilde{x} = \left(x_1,\cdots,x_m\right)$, $\tilde{U}= \left\{\left(x_1,\cdots,x_m\right) : x\in U\right\}$ and $f : \tilde{U}\ra \mathbb{R}^{n-m}$ is $C^p$. We may assume also that there exists $C > 0$ such that $\norm Df\left(\tilde{x}\right)\norm\leq C$ for every $\tilde{x}\in\tilde{U}$.
		
		Let $x, y \in M\cap U$ and consider the curve $\gamma$ in $M$ defined by $\left\{\left(\tilde{z}, f\left(\tilde{z}\right)\right) : \tilde{z}\in [\tilde{x},\tilde{y}]\right\}$, where $[\tilde{x},\tilde{y}]$ is the segment joining $\tilde{x}$, $\tilde{y}$ in $\tilde{U}\subset \mathbb{R}^m$. Hence we have that
		\[
		d(x,y)\leq \mathcal{L}(\gamma) \leq K\norm \tilde{x}-\tilde{y}\norm\leq K\norm x-y\norm,
		\]
		since  $\norm \dot{\gamma}(t)\norm\leq K:=\sqrt{1+C^2}$ for all $t\in \left[0,\norm \tilde{x}-\tilde{y}\norm\right]$.
	\end{proof}
	\begin{theorem}\label{thm4.1}
		Let $S$ be a closed subset of $M$ and $x\in S$. Then $\mathbf{N}^P_S(x)=N^P_S(x)\oplus N_x$.
	\end{theorem}
	\begin{proof}
		Let us denote $\mathbf{v} := v + n$  where $v\in N^P_S(x)$ and $n\in N_x$. Then there exist positive constants $\eta$ and $\delta$ such that
		\[
		\li v, \exp_x^{-1}y \ri \leq\delta\;d^2(x,y)=\delta\norm\exp_x^{-1}y\norm^2,
		\]
		for every $y\in B(x,\eta)\cap S$. Hence for every $y\in \mathbf{B}(x,\eta)\cap S$ we have
		\begin{align*}
			\li \mathbf{v}, y-x \ri = & \li \mathbf{v}, \exp_x^{-1}y\ri + \li \mathbf{v}, y- \left(x + \exp_x^{-1}y\right)\ri \\
			= & \li v, \exp_x^{-1}y\ri + \li \mathbf{v}, y-\left(x + \exp_x^{-1}y\right)\ri\\
			\leq & \delta\norm\exp_x^{-1}y\norm^2+\norm \mathbf{v}\norm\left\norm y-\exp_x^{-1}y- \left(x-\exp_x^{-1}x\right)\right \norm
		\end{align*}
		We now obtain an upper bound for the expression $\norm y-\exp_x^{-1}(y)-(x-\exp_x^{-1}(x))\norm$ using Taylor's expansion.
		To this end, we consider the map $f: B(x,\eta)\ra \mathbb{R}^n$ defined by $f(y):=I(y)-\exp_x^{-1}(y)$, where $B(x,\eta)\subset M$ is a convex ball.
		Let $y\in B(x,\eta)$ and
		\[
		\gamma(t):=\exp_x\left(t\frac{\exp_x^{-1}(y)}{\norm \exp_x^{-1}(y)\norm}\right) \quad \forall\  t\in [0,d(x,y)],
		\]
		be the geodesic with unit speed joining $x,y$. We put $g(t):=f\left(\gamma(t)\right)$ for all $t\in [0,d(x,y)]$. Then we have $g'(0)=0$ and
		\begin{align*}
			g''(t)= & \frac{d^2}{dt^2}\left(I(\gamma(t))-\exp_x^{-1}(\gamma(t))\right) \\
			= & \frac{d}{dt}\left(\dot{\gamma}(t)-D\exp_x^{-1}(\gamma(t))(\dot{\gamma}(t))\right) \\
			= & \ddot{\gamma}(t)-D^2\exp_x^{-1}(\gamma(t))(\dot{\gamma}(t))^2.
		\end{align*}
		Note that since $\gamma$ is a geodesic on $M$,  $\ddot{\gamma}(t)\in N_{\gamma(t)}$ and  we have $\ddot{\gamma}=\two\left(\dot{\gamma},\dot{\gamma}\right)$ where $\two$ is the second fundamental form of $M$. Then we have
		\[
		g''(t_0)=\two_z\left(w,w\right)-D^2\exp_x^{-1}(z)(w)^2,
		\]
		where $z:=\gamma(t_0)$, $w:=\dot{\gamma}(t_0)$ and $\two_z:T_zM\times T_zM\ra (T_zM)^{\bot}$ is a bilinear form.
		
		Therefore there exist $z\in B(x,\eta)$ and $w\in T_zM$ with $\norm w\norm=1$ such that
		\[
		f(y)=f(x)+\frac{1}{2}\left(\two_z\left(w,w\right)-D^2\exp_x^{-1}(z)(w)^2\right)d^2(x,y).
		\]
		It follows that
		\begin{align*}
			\norm y-\exp_x^{-1}(y)-(x-\exp_x^{-1}(x))\norm\leq & \frac{1}{2}\left\norm\two_z-D^2\exp_x^{-1}(z)\right\norm d^2(x,y)\\
			\leq & \frac{1}{2}\left(\left\norm\two_z\right\norm+\left\norm D^2\exp_x^{-1}(z)\right\norm \right)d^2(x,y).
		\end{align*}
		Hence  we conclude that
		\[
		\li \textbf{v},y-x\ri\leq \left(\delta+\frac{\norm \textbf{v} \norm}{2}\left(\kappa+C\right)\right)d^2(x,y),
		\]
		where $C:= \sup_{z\in B(x,\eta)}\norm D^2\exp_x^{-1}(z)\norm$ and $\kappa$ is a positive constant such that  $\left\norm\two_z\right\norm\leq \kappa$ for all $z\in B(x,\eta)$. Finally, using Lemma \ref{lem4.1} we deduce that
		\begin{equation}\label{eq4.1}
			\li \textbf{v},y-x\ri\leq K\left(\delta+\frac{\norm \textbf{v} \norm}{2}\left(\kappa+C\right)\right)\norm x-y\norm^2,
		\end{equation}
		for every $y\in \mathbf{B}(x,\eta)\cap S$. It follows that $N^P_S(x)\oplus N_x\subseteq \mathbf{N}^P_S(x)$.
		
		Conversely, if $\textbf{v}\in \mathbf{N}^P_S(x)$ , then there exist positive constants $\eta$ and $\delta$ such that
		\[
		\li \textbf{v},z-x\ri\leq \delta \norm z-x\norm^2,
		\]
		for every $z\in \mathbf{B}(x,\eta)\cap S$. Moreover, $\textbf{v}=v+n$ where $v\in T_xM$ and $n\in N_x$. In particular, for every $y\in S$ with the property that $\norm \exp_x^{-1}(y)\norm <\eta$, putting $z=x+\exp_x^{-1}(y)$, we obtain
		\[
		\li v, \exp_x^{-1}y \ri \leq\delta \norm \exp_x^{-1}(y)\norm^2,
		\]
		which implies that $v\in N^P_S(x)$.
	\end{proof}
	\begin{theorem}\label{thm4.2}
		If $S \subset M$ is prox-regular, then it is prox-regular as a subset of $\mathbb{R}^n$.
	\end{theorem}
	\begin{proof}
		It is an immediate consequence of Theorem \ref{thm4.1}, since all
		the constants in the inequality \eqref{eq4.1} depend continuously on $x$.
	\end{proof}
	\begin{proposition}
		If $\gamma:I\ra M$ is a weak geodesic on $S$ as a subset of $M$, then it is a weak geodesic
		in $S$ as a subset of $\mathbb{R}^n$.
	\end{proposition}
	\begin{proof}
		It is a consequence of the fact that $D_t\dot{\gamma}(t) =P_{T_{\gamma(t)}M} \left(\ddot{\gamma}(t)\right)$ if we look $S$ as a
		subset of $\mathbb{R}^n$. Then Theorem \ref{thm4.1} gives us the result.
	\end{proof}
	Hence according to  \cite[Theorem 2.7]{Canino1} we get the following result.
	\begin{corollary}\label{coro4.1}
		If $\gamma:I\ra M$ is a weak geodesic on $S$, then $\gamma$ is locally minimizing in $S$.
	\end{corollary}
	
	\begin{example}
		Let $M$ be the 2-dimensional unit sphere, $x^2+y^2+z^2=1$ in $\mathbb{R}^3$ with the induced metric
		$g^\circ=d\theta ^2+\sin ^2 \theta\:d\phi ^2$ in spherical coordinates $\left(\theta,\phi\right)$ defined on
		$M\setminus\left\{(x,y,z): x\leq 0,y=0\right\}$ by
		\[
		(x,y,z)=\left(\sin \theta \cos \phi, \sin \theta\sin \phi,\cos \theta\right) \quad 0<\theta<\pi,-\pi<\phi<\pi.
		\]
		
		Let $\theta_0$ be such that $\frac{\pi}{2}<\theta_0<\pi$ and
		$S_1$, $S_2$ be the closed subsets of $M$ defined by
		\[
		S_1:=\left\{\left(\theta,\phi\right): 0\leq \theta \leq \theta_0\right\},\ \ \
		S_2:=\left\{\left(\theta,\phi\right): \theta_0\leq \theta \leq\pi\right\}.
		\]
		By \cite[Theorem 4.18]{convex}, $S_1$ is a prox-regular subset of $M$. On the other hand, $S_2$ is a convex set  in $M$ and
		hence it is prox-regular.
		Using \cite[Theorem 1]{Minimizing}, the proximal normal cones to $S_1$ and $S_2$ at each point $x\in \partial S_1=\partial S_2$
		are obtained as
		\[
		N^P_{S_1}(x)=\cone \{(1,0)\}=\{\lambda\:\partial/\partial \theta:\lambda\geq 0\},
		\]
		and
		\[
		N^P_{S_2}(x)=\cone \{(-1,0)\}=\{-\lambda\:\partial/\partial \theta:\lambda\geq 0\}.
		\]
		
		We now consider the unit speed curve $\gamma$ on $M$ defined by
		\[
		\gamma(t):=\left(\theta_0,\frac{t}{\sin \theta_0}\right) \qquad \forall t\in I:=\left[-\pi/2\sin \theta_0,\pi/2\sin \theta_0\right],
		\]
		and we have
		\[
		D_t\dot{\gamma}\left(t\right)=\left(-\frac{\cos \theta_0}{\sin \theta_0}\right)\frac{\partial}{\partial \theta}
		\qquad \forall t\in I.
		\]
		Then $\gamma$ is a weak geodesic on $S_1$ and hence it is locally minimizing in $S_1$, but $\gamma$ is neither
		locally minimizing in $S_2$ nor a weak geodesic on $S_2$.
	\end{example}

	We now assume that $M$ is an m-dimensional complete Riemannian manifold and $S$ is a nonempty closed subset of $M$. We prove that our desired results are
	invariant under $C^{\infty}$ isometries.  By the Nash embedding theorem \cite{Nash}, let
	$\phi : M \ra  \mathbb{R}^n$ be an isometric embedding and $\tilde{M}:=\phi(M)$. Then $\tilde{M}$
	is a complete Riemannian submanifold of $\mathbb{R}^n$ and $\phi(S)$ is a nonempty closed subset of $\tilde{M}$.
	
	\begin{lemma}\label{lem4.2}
		If $x\in S$, then
		\[
		N^P_{\phi(S)}\left(\phi(x)\right)= D\phi(x)\left(N^P_S(x)\right).
		\]
	\end{lemma}
	\begin{proof}
		Let $v\in N^P_S(x)$.
		Since $\phi$   transforms locally geodesics passing through $x$ on $M$ into geodesics passing through $\phi(x)$ on $\tilde{M}$,
		there exists a positive constant $\delta$ such that
		\begin{equation}\label{ineqdel}
			\li D\phi(x)(v), \exp_{\phi(x)}^{-1}\left(\phi(y)\right)\ri=  \li D\phi(x)(v), D\phi(x)\left(\exp_{x}^{-1}y\right)\ri
		\end{equation}
		\[
		=  \li v, \exp_x^{-1}y\ri \leq \delta d_M^2(x,y)=\delta d_{\tilde{M}}^2\left(\phi(x),\phi(y)\right),
		\]
		for every $y\in U\cap S$, where $U$ is a convex neighborhood of $x$ in $M$.
		
	\end{proof}
	
	This result allows us to prove the following theorems.
	
	\begin{theorem}\label{thm4.3}
		If $S \subseteq M$ is prox-regular, then $\phi(S)\subseteq \tilde{M}$ is prox-regular.
	\end{theorem}
	\begin{proof}
			Using the inequality \eqref{ineqdel} and by the choice of $\delta$, we get the desired result.
	\end{proof}
	\begin{theorem}\label{thm4.4}
		If $\gamma$ is a weak geodesic on $S$, then $\phi\circ\gamma$ is a weak geodesic on $\phi(S)$ as a subset of $\tilde{M}$.
	\end{theorem}
	\begin{proof}
		We  show that
		\[
		\tilde{D}_t\left(\dot{\phi\circ\gamma}\right)(t)\in N^P_{\phi(S)}\left(\left(\phi\circ\gamma\right)(t)\right) \qquad a.\ e.\ t\in I,
		\]
		where $\tilde{D}_t$ denotes the covariant derivative on $\tilde{M}$. By Lemma \ref{lem4.2}, this is equivalent to
		\[
		D\phi\left(\gamma(t)\right)^{-1}\left[\tilde{D}_t\left(\dot{\phi\circ\gamma}\right)(t) \right]\in N^P_S\left(\gamma(t)\right) \qquad a.\ e.\ t\in I.
		\]
		Let us calculate $\tilde{D}_t\left(\dot{\phi\circ\gamma}\right)(t)$. It can be seen as a vector in the Euclidean
		space where $\tilde{M}$ lies, that is the orthogonal projection of $\left(\phi\circ\gamma\right)''(t)$ onto $T_{\left(\phi\circ\gamma\right)(t)}\tilde{M}$. On the other hand, we have
		\begin{align*}
			\left(\phi\circ\gamma\right)''(t)= & \left(D\phi\left(\gamma(t)\right)\left(\dot{\gamma}(t)\right) \right)' \\
			= & D^2\phi\left(\gamma(t)\right)\left(\dot{\gamma}(t)\right)^2+D\phi\left(\gamma(t)\right)\left(D_t\dot{\gamma}(t)\right),
		\end{align*}
		where $D\phi\left(\gamma(t)\right)\left(D_t\dot{\gamma}(t)\right)\in T_{\left(\phi\circ\gamma\right)(t)}\tilde{M}$.
		
		We now calculate $D\phi\left(\gamma(t)\right)^{-1}\left[\tilde{D}_t\left(\dot{\phi\circ\gamma}\right)(t) \right]$. The second term is $D_t\dot{\gamma}(t)$  which belongs to $N^P_S\left(\gamma(t)\right)$. For the first term, we observe that if $\sigma$ is a geodesic on
		$M$ passing through $\gamma(t)$ such that $\dot{\sigma}(0) = \dot{\gamma}(t)$, then $\phi\circ\sigma$
		is a geodesic on $\tilde{M}$ passing through $\phi\left(\gamma(t)\right)$. This implies that $\left(\phi\circ\sigma\right)''(0)$
		is normal to $\tilde{M}$ at $\phi\left(\gamma(t)\right)$. On the other hand, we have
		\begin{align*}
			\left(\phi\circ\sigma\right)''(0)= & D^2\phi\left(\gamma(t)\right)\left(\dot{\gamma}(t)\right)^2+D\phi\left(\gamma(t)\right)\left(D_t\dot{\sigma}(0)\right) \\
			= & D^2\phi\left(\gamma(t)\right)\left(\dot{\gamma}(t)\right)^2.
		\end{align*}
		Then we conclude that
		\[
		D^2\phi\left(\gamma(t)\right)\left(\dot{\gamma}(t)\right)^2\in N_{\left(\phi\circ\gamma\right)(t)}\tilde{M},
		\]
		and hence $\tilde{D}_t\left(\dot{\phi\circ\gamma}\right)(t)=D\phi\left(\gamma(t)\right)\left(D_t\dot{\gamma}(t)\right)$.
		It follows that
		\[
		D\phi\left(\gamma(t)\right)^{-1}\left[\tilde{D}_t\left(\dot{\phi\circ\gamma}\right)(t) \right]=D_t\dot{\gamma}(t)\in N^P_S\left(\gamma(t)\right).
		\]
	\end{proof}
	
	Therefore using Corollary \ref{coro4.1}, we obtain the following result.
	\begin{corollary}
		Let $M$ be a complete Riemannian manifold and $\gamma:I\ra M$ be a weak geodesic on $S$. Then $\gamma$ is locally minimizing in $S$.
	\end{corollary}
	
	\section{Closed weak geodesics on prox-regular sets}\label{sec5}
	In this section, we define closed weak geodesics on prox-regular subsets of Riemannian manifolds and we prove that these are the only viscosity critical points of the energy functional; \re{see \cite{weak} for the discussion on weak geodesics.} Finally, we discuss the local minimality properties of closed weak geodesics.
	
	Throughout  this section, we suppose that $S$ is a nonempty and closed  $\varphi$-convex subset of $M$ and $U$ is an open neighborhood of $S$ on which $P_S$ is single-valued and locally Lipschitz.
	Motivated by  \cite{Canino,weak}, a curve $\gamma\in H^1(I,M)$ with $\im(\gamma)\subseteq S$ is called a closed weak geodesic on $S$ if $\gamma$ is a weak geodesic with the properties that $\gamma(0)=\gamma(1)$ and  $\dot{\gamma}_{+}(0)=\dot{\gamma}_{-}(1)$. More precisely, we give the following definition.
	\begin{definition}
		 A closed weak geodesic on $S$ is a continuous curve $\gamma : I\ra M$  which satisfies the following properties:
		\begin{itemize}
			\item [(i)] $\gamma(t)\in S \quad \forall t\in I$;
			\item [(ii)] $\gamma\in W^{2,2}(I,M)$;
			\item [(iii)] $D_t\dot{\gamma}(t)\in N^P_S\left(\gamma(t)\right) \quad a.\ e.\ t\in I$;
			\item [(iv)] $\gamma(0)=\gamma(1)$ and  $\dot{\gamma}_{+}(0)=\dot{\gamma}_{-}(1)$.
		\end{itemize}
	\end{definition}
	
	We now prove that closed weak geodesics on $S$ are the only nonsmooth critical  points for the appropriately defined energy functional. For this purpose, we consider the following set of admissible curves:
	\[
	\mathcal{C}:=\left\{\eta\in H^1(I,M) : \eta(t)\in S, \ \forall t\in I, \eta(0)= \eta(1)\right\}.
	\]
	In the special case, when the set $S$ is all of $M$, the set $\mathcal{C}$ is a complete submanifold of $H^1(I,M)$, see \cite{Klingenberg}. We now consider the energy functional as follows:
	\[
	\hspace{-3.8cm} f_{{}_c}: L^2(I,M)\ra \mathbb{R}\cup \{+\infty\}
	\]
	
	\[
	f_{{}_c}(\eta)=\left\{
	\begin{array}{lr}
		\frac{1}{2} \int_0^1 \norm \dot{\eta}(t)\norm ^2 dt & \eta\in \mathcal{C} \\
		& \\
		+\infty & \eta\in L^2(I,M)\setminus\mathcal{C}.
	\end{array}\right.
	\]
	
	\begin{proposition}
		The energy functional $f_{{}_c}$ is lower semicontinuous and  $\mathcal{C}\subseteq L^2(I,M)$ is closed.
	\end{proposition}
	\begin{proof}
		We consider the map $g: \mathcal{D}\ra \mathbb{R}\times\mathbb{R}$ defined by
		\[
		g(\eta):=\left(d_{\infty}\left(\eta,P_S\left(\eta\right)\right),d\left(\eta(0),\eta(1)\right)\right),
		\]
		where
		\[
		\mathcal{D}:=\left\{\eta\in H^1(I,M) : \eta(t)\in U, \ \forall t\in I\right\}
		\]
		is open in $H^1(I,M)$ and $P_S:\mathcal{D}\ra \mathcal{D}$ defined by $P_S\left(\eta\right)(t):=P_S\left(\eta(t)\right)$ for all $t\in I$ is a continuous functional map; see \cite[Proposition 3.3]{weak}. Then $g$ is continuous and hence $\mathcal{C}=g^{-1}\{(0,0)\}$ is closed.
		
		Therefore $f_{{}_c}=E\times 1_{\mathcal{C}}$ is lower semicontinuous where $1_{\mathcal{C}}$ is the characteristic function  of $\mathcal{C}$ and $E:H^1(I,M)\ra \mathbb{R}$ defined by
		\[
		E\left(\eta\right):=\frac{1}{2}\int_{0}^{1}\norm \dot{\eta}(t)\norm^2dt,
		\]
		is smooth; see \cite{Klingenberg}.
	\end{proof}
	\begin{lemma}\label{lem1}
		Let $\gamma\in \mathcal{C}$ and $\xi\in D^-f_{{}_c}(\gamma)$. Then
		\[
		\int_{0}^{1}\li \dot{\gamma}, D_tV\ri dt\geq \int_{0}^{1}\li \xi, P_{\gamma}V\ri dt-\int_{0}^{1}\left(2\varphi(\gamma)+\tau\right)\norm V-P_{\gamma}V\norm \norm \dot{\gamma}\norm^2 dt,
		\]
 for every vector field  $V\in H^1\left(I, \gamma^{-1}TM\right)$ with the property that $V(0)=V(1)$, where $\tau=\tau\left(\gamma\right)$ is a piecewise constant function  on $I$.
	\end{lemma}
	\begin{proof}
		Let $V\in H^1\left(I, \gamma^{-1}TM\right)$ and $V(0)=V(1)$. Using the vector field $V$, we construct some variations $\Gamma$ and $\tilde{\Gamma}$ of $\gamma$ as follows
		\[
		\Gamma_s(t):=\exp_{\gamma(t)}\left(sV(t)\right) \quad \forall t\in I,
		\]
		and
		\[
		\tilde{\Gamma}_s(t):= P_S\left(\Gamma_s(t)\right) \quad \forall t\in I,
		\]
		for  sufficiently small $s$. Hence $\Gamma_s, \tilde{\Gamma}_s\in H^1(I,M)$. Moreover, we have $\tilde{\Gamma}_s(0)=\tilde{\Gamma}_s(1)$ and hence $\tilde{\Gamma}_s\in \mathcal{C}$.
		Therefore 
		we obtain that
		\[
		\hspace{-6cm}\int_{0}^{1}\li \dot{\gamma}, D_tV\ri dt  -\int_{0}^{1}\li \xi, P_{\gamma}V\ri dt\,=
		\]
		\[
		\lim_{s\ra 0^{+}}\frac{1}{s}\int_{0}^{1}\left(\frac{1}{2}\norm \dot{\Gamma_s}(t)\norm^2-\frac{1}{2}\norm \dot{\gamma}(t)\norm^2-\left\li \xi , \exp_{\gamma(t)}^{-1}P_S\left(\exp_{\gamma(t)}sV(t)\right)\right\ri\right)dt
		\]
		\[
		\geq \liminf_{s\ra 0^{+}}\frac{1}{s}\int_{0}^{1}\left(\frac{1}{2}\norm \dot{\tilde{\Gamma}}_s(t)\norm^2-\frac{1}{2}\norm \dot{\gamma}(t)\norm^2-\left\li \xi , \exp_{\gamma(t)}^{-1}\tilde{\Gamma}_s(t)\right\ri\right)dt
		\]
		\[
		+\liminf_{s\ra 0^{+}}\frac{1}{2s}\int_{0}^{1}\left(\norm \dot{\Gamma_s}(t)\norm^2-\norm \dot{\tilde{\Gamma}}_s(t)\norm^2\right)dt.
		\]
		Since $\xi\in D^-f_{{}_c}(\gamma)$, the first term is nonnegative and the second term has been estimated from below in \cite[Theorem 3.7]{weak}. Then we conclude that
		\[
		\int_{0}^{1}\li \dot{\gamma}, D_tV\ri dt- \int_{0}^{1}\li \xi, P_{\gamma}V\ri dt\geq -\int_{0}^{1}\left(2\varphi(\gamma)+\tau\right)\norm V-P_{\gamma}V\norm \norm \dot{\gamma}\norm^2 dt.
		\]
	\end{proof}
	
	Note that a special case of Lemma \ref{lem1}, where $\gamma$ is a closed curve on a $\varphi$-convex subset of a Hilbert space, has been proved in \cite{Canino} with $\tau\equiv 0$.

	The following theorem gives some properties of closed curves containing in the domain of $D^-f_{{}_c}$ and states the relationship between the subderivatives of $f_{{}_c}$ at $\gamma\in \mathcal{C}$ and the proximal normal vectors of $S$ at the image of $\gamma$.
	\begin{theorem}\label{thm1}
		Let $\gamma\in \mathcal{C}$. Then $D^-f_{{}_c}(\gamma)\neq \emptyset$ if and only if $\gamma\in W^{2,2}(I,M)$ and $\dot{\gamma}_{+}(0)=\dot{\gamma}_{-}(1)$. Moreover, if $\xi\in L^2\left(I,\gamma^{-1}TM\right)$, then $\xi\in D^-f_{{}_c}(\gamma)$ if and only if
		\[
		\xi(t)+D_t\dot{\gamma}(t)\in N^P_S\left(\gamma(t)\right),\ a.\ e. \ t\in I.
		\]
	\end{theorem}
	\begin{proof}
		Let $\xi\in D^-f_{{}_c}(\gamma)$.  Using Lemma \ref{lem1} we have
		\begin{equation}\label{eq1}
			\int_{0}^{1}\li \dot{\gamma}, D_tV\ri dt\geq \int_{0}^{1}\li \xi, P_{\gamma}V\ri dt-\int_{0}^{1}\left(2\varphi(\gamma)+\tau\right)\norm V-P_{\gamma}V\norm \norm \dot{\gamma}\norm^2 dt,
		\end{equation}
		for all  $V\in H^1\left(I, \gamma^{-1}TM\right)$ with $V(0)=V(1)$  where $\tau$ is a piecewise constant function  on $I$.
		This inequality implies that
		\[
		\norm D_t\dot{\gamma}\norm_{L^2}\leq \left(1+\left(2\bar{\varphi}+C\right)\norm \dot{\gamma}\norm_{L^2}\right)\left(\norm \xi\norm_{L^2}+\left(2\bar{\varphi}+C\right)\norm \dot{\gamma}\norm_{L^2}^2\right),
		\]
		and hence $\gamma\in W^{2,2}(I,M)$.
		
		We now prove that  $\dot{\gamma}_{+}(0)=\dot{\gamma}_{-}(1)$. It suffices to show that
		\[
		\left\li  \dot{\gamma}_{-}(1)-\dot{\gamma}_{+}(0), v\right\ri\geq 0 \qquad\forall v\in T_{\gamma(0)}M.
		\]
		For given $v\in T_{\gamma(0)}M$, we find an $H^1$-vector field $V$ along $\gamma$ with the property that  $V(0)=V(1)=v$.
		Indeed, this vector field can be constructed as follows.
		Let $\left(U_{\alpha}, x^1,\cdots,x^n\right)$ be a chart of $M$ around $\gamma(0)$ and $v=\sum_{i=1}^{n}v_i\frac{\partial}{\partial x^i}|_{\gamma(0)}$, for some $v_i\in \mathbb{R}$, $i=1,\cdots,n$. We choose $t_0\in I$ with the properties that $t_0<1$ and $\gamma(t)\in U_{\alpha}$ for all $t\in [t_0,1]$. Putting $w:=L^{\gamma}_{\gamma(0),\gamma(t_0)}v$, we have $w\in T_{\gamma(t_0)}M$ and let $w=\sum_{i=1}^{n}w_i\frac{\partial}{\partial x^i}|_{\gamma(t_0)}$ for some $w_i\in \mathbb{R}$, $i=1,\cdots,n$. We now define
		\[
		V(t):=\left\{
		\begin{array}{lr}
			L^{\gamma}_{\gamma(0),\gamma(t)}v & 0\leq t\leq t_0 \\
			& \\
			\sum_{i=1}^{n}\frac{1}{t_0-1}\left( (t-1)w_i+(t_0-t)v_i\right)\frac{\partial}{\partial x^i}|_{\gamma(t)} & t_0\leq t\leq 1.
		\end{array}
		\right.
		\]
		Clearly, $V\in H^1\left(I, \gamma^{-1}TM\right)$ and $V(0)=V(1)=v$.
		
		We now consider a sequence of $H^1$-vector fields $V_n$ along $\gamma$ defined by $V_n(t):=u_n(t)V(t)$ for all $t\in I$ where the functions $u_n\in W^{1,2}(0,1)$ have the following properties
		\[
		0\leq u_n\leq 1,\quad u_n(0)=u_n(1)=1,\  \text{and}\ \ u_n(t)=0 \quad \forall t\in [\frac{1}{2n},1-\frac{1}{2n}].
		\]
		Then $V_n(0)=V_n(1)$ and putting $V_n$ in \eqref{eq1}, we have
		\[
		\hspace{-2cm}\int_{0}^{1}u'_n\li \dot{\gamma}, V\ri dt+\int_{0}^{1}u_n\li \dot{\gamma}, D_tV\ri dt\geq \int_{0}^{1}u_n\li \xi, P_{\gamma}V\ri dt
		\]
		\[
		\hspace{4cm}-\int_{0}^{1}u_n\left(2\varphi(\gamma)+\tau\right)\norm V-P_{\gamma}V\norm \norm \dot{\gamma}\norm^2 dt.
		\]
		Integrating by parts in the first term implies that
		\[
		\hspace{-1cm}\left\li  \dot{\gamma}_{-}(1)-\dot{\gamma}_{+}(0), v\right\ri-\int_{0}^{1}u_n\frac{d}{dt}\li \dot{\gamma}, V\ri dt+\int_{0}^{1}u_n\li \dot{\gamma}, D_tV\ri dt\geq
		\]
		\[
		\hspace{1cm}\int_{0}^{1}u_n\li \xi, P_{\gamma}V\ri dt-\int_{0}^{1}u_n\left(2\varphi(\gamma)+\tau\right)\norm V-P_{\gamma}V\norm \norm \dot{\gamma}\norm^2 dt.
		\]
		Then passing to the limit as $n\ra \infty$, we conclude that
		\[
		\left\li  \dot{\gamma}_{-}(1)-\dot{\gamma}_{+}(0), v\right\ri\geq 0.
		\]
		
		In order to prove
		\[
		D_t\dot{\gamma}(t)+\xi(t)\in N^P_S\left(\gamma(t)\right)\quad  \ a.\ e.\ t\in I,
		\]
		we assume that $t_0\in I$ and $\omega\in T^B_S\left(\gamma(t_0)\right)$. 
		We now consider a sequence of vector fields $V_n\in H^1_0\left(I,\gamma^{-1} TM\right)$ along $\gamma$ defined by
		\[
		V_n(t):=u_n\left(t-t_0\right)L_{t_0,t}(\omega) \quad \forall\ t\in I,
		\]
		where  $u_n\in C^{\infty}_0\left(\mathbb{R}\right)$ is a sequence of functions with the properties that
		\[
		u_n\geq 0,\quad \supp u_n\subseteq \left[-\frac{1}{n},\frac{1}{n}\right], \quad \int u_n=1,\quad \forall n\in \mathbb{N}.
		\]
		Then using Lemma \ref{lem1}, we obtain that
		\[
		-\left\li\int_{0}^{1} u_n\left(t-t_0\right)L_{t,t_0}\left(D_t\dot{\gamma}\right)dt, \omega\right\ri\geq \int_{0}^{1}u_n\left(t-t_0\right)\li \xi, P_{\gamma}\left(L_{t_0,t}\:\omega\right)\ri dt
		\]
		\[
		-\left(2\bar{\varphi}+C\right)\norm \dot{\gamma}\norm^2_{\infty}\int_{0}^{1}u_n\left(t-t_0\right)\norm L_{t_0,t}\:\omega-P_{\gamma}\left(L_{t_0,t}\:\omega\right)\norm  dt.
		\]
		As $n\ra \infty$, this implies  that
		\[
		\left\li D_t\dot{\gamma}(t_0)+\xi(t_0),\omega\right\ri \leq 0,
		\]
		and then we get the result.
		
		To prove the converse statement, let $\gamma\in \mathcal{C}\cap W^{2,2}(I,M)$,  $\dot{\gamma}_{+}(0)=\dot{\gamma}_{-}(1)$ and $\xi\in L^2\left(I,\gamma^{-1}TM\right)$ be such that
		\[
		\xi+D_t\dot{\gamma}\in N^P_S\left(\gamma\right),\ a.\ e.
		\]
		We Suppose that $W$ is an open neighborhood  of $\im \left(\gamma\right)$ with compact closure and $\delta\leq\se\leq\Delta$ on $W$ with $\Delta\geq 0$. Then for every $\eta\in \mathcal{C}$ with $d_{\infty}\left(\gamma, \eta\right)< \bar{r}$ we have
		\begin{equation}\label{ineq5}
			\frac{1}{2}\int_{0}^{1}\norm \dot{\eta}(t)\norm^2dt\geq  \frac{1}{2}\int_{0}^{1}\left(c\left(\gamma,\eta\right)-1\right)\norm \dot{\gamma}(t)\norm^2dt+\int_{0}^{1}\left\li\xi,\exp_{\gamma}^{-1}\eta\right\ri dt
		\end{equation}
		\[
		-d_{\infty}(\gamma ,\eta )\norm\exp_{\gamma}^{-1}\eta\norm_{L^{2}}
		\big( \bar{\varphi}\norm\xi+
		D_t\dot{\gamma}\norm_{L^{2}}+\frac{1}{2}|R|_{\infty}\norm \dot{\gamma}\norm_{L^{\infty}}\norm\dot{\eta}\norm_{L^{2}}\big),
		\]
		where $c\left(\gamma,\eta\right)=2\sqrt{\Delta}\:d\left(\gamma,\eta\right)
		\cot\left(\sqrt{\Delta}\:d\left(\gamma,\eta\right)\right)$, $\bar{r}=\min_{t\in I}r\left(\gamma(t)\right)$ and $R$ is the curvature tensor on $M$. Indeed, let  $\eta\in \mathcal{C}$ and $d_{\infty}\left(\gamma, \eta\right)< \bar{r}$, hence 
		we have
		\[
		\frac{1}{2}\norm \dot{\eta}\norm^2-\frac{c\left(\gamma,\eta\right)-1}{2}\norm \dot{\gamma}\norm^2\geq \frac{1}{2}\norm L_{\eta,\gamma}\dot{\eta}-\dot{\gamma}\norm^2+\frac{d}{dt}\left\li \dot{\gamma}, \exp_{\gamma}^{-1}\eta\right\ri
		\]
		\[
		-  \left\li D_t\dot{\gamma}, \exp_{\gamma}^{-1}\eta\right\ri-\frac{1}{2}|R|_{\infty} d^2\left(\gamma, \eta\right)\norm \dot{\gamma}\norm \norm \dot{\eta}\norm.
		\]
		Then integrating from both side and using the assumptions $\xi+D_t\dot{\gamma}\in N^P_S\left(\gamma\right),\ a.\ e.$ and $\dot{\gamma}_{+}(0)=\dot{\gamma}_{-}(1)$, we conclude that
		\[
		\frac{1}{2}\int_{0}^{1}\norm \dot{\eta}(t)\norm^2dt -  \frac{1}{2}\int_{0}^{1}\left(c\left(\gamma,\eta\right)-1\right)\norm \dot{\gamma}(t)\norm^2dt-\int_{0}^{1}\left\li\xi,\exp_{\gamma}^{-1}\eta\right\ri dt
		\]
		\[
		\geq \frac{1}{2}\int_{0}^{1}\norm L_{\eta,\gamma}\dot{\eta}-\dot{\gamma}\norm^2dt-\int_{0}^{1}\varphi\left(\gamma\right)\norm\xi+
		D_t\dot{\gamma}\norm\norm\exp_{\gamma}^{-1}\eta\norm^2 dt
		\]
		\[
		-\frac{1}{2}|R|_{\infty}\int_{0}^{1} \norm \dot{\gamma}\norm \norm \dot{\eta}\norm\norm\exp_{\gamma}^{-1}\eta\norm^2 dt.
		\]
		Thus the inequality \eqref{ineq5} is obtained; see  \cite{weak}.
		Moreover, we have
		\[
		\Big|  \int_{0}^{1}\left(c\left(\gamma,\eta\right)-2\right)\norm \dot{\gamma}(t)\norm^2dt\Big| \leq
		K\left\norm \exp_{\gamma}^{-1}\eta\right\norm_{L^2}^2
		\]
		for a suitable constant $K$ and this implies that
		\[
		\liminf_{d_{\infty}\left( \eta,\gamma\right)\ra 0}\frac{f_{{}_c}\left(\eta\right)-f_{{}_c}\left(\gamma\right)-\left\li\xi,\exp_{\gamma}^{-1}\eta\right\ri_{L^2} }{\left\norm \exp_{\gamma}^{-1}\eta\right\norm_{L^2}}\geq 0.
		\]
		It follows that $\xi\in D^-f_{{}_c}(\gamma)$.

		In particular, since
		\[
		-D_t\dot{\gamma}+D_t\dot{\gamma}\in N^P_S\left(\gamma\right),\ a.\ e.,
		\]
		we have $-D_t\dot{\gamma}\in D^-f_{{}_c}(\gamma)$ and this completes the proof.
	\end{proof}

	\begin{theorem}
		Let $\gamma\in \mathcal{C}$, then $0\in D^-f_{{}_c}(\gamma)$ if and only if $\gamma$ is a closed weak geodesic on $S$.
	\end{theorem}
	\begin{proof}
		Using Theorem \ref{thm1}, we have $0\in D^-f(\gamma)$ if and only if $\gamma\in W^{2,2}(I,M)$, $\dot{\gamma}_{+}(0)=\dot{\gamma}_{-}(1)$ and $D_t\dot{\gamma}(t)\in N^P_S\left(\gamma(t)\right)$ for almost all $t\in I$.
	\end{proof}
	
	Finally, in the following theorem we study the local minimality of closed weak geodesics on $S$.
	
	
	
	\begin{theorem}
		Let $\gamma\in \mathcal{C}$.  If $\gamma$ is a locally minimizing curve in $S$ parameterized with constant speed and there exists $h>0$ such that $f(\gamma)< f(\eta)$ for every  $\eta\in \mathcal{C}\setminus \{\gamma\}$ with the property that $\eta=\gamma$ on $ \left[h,1-h\right]$, then $\gamma$   is a closed weak geodesic on $S$.
		
		The converse statement holds if  the sectional curvatures of $M$ are nonpositive in a neighborhood of $\im(\gamma)$.
	\end{theorem}
	\begin{proof}
		Let $\gamma$ be a locally minimizing curve with constant speed in $S$. Then according to Theorem \ref{localweak}, the curve $\gamma$   is a weak geodesic on $S$. We now show that $\dot{\gamma}_{+}(0)=\dot{\gamma}_{-}(1)$.  To this end, 
		for fixed $t_0\in (0,1)$ we define a curve $\alpha\in \mathcal{C}$ defined by
		\[
		\alpha(t):=\left\{
		\begin{array}{ll}
			\gamma\left(t+t_0\right) & 0\leq t\leq 1-t_0  \\
			\gamma\left(t-1+t_0\right) & 1-t_0\leq t\leq 1.
		\end{array}
		\right.
		\]
		Then by the assumption, $\alpha$ is a locally minimizing curve with constant speed in $S$ and hence it is a weak geodesic on $S$. It follows that $\alpha\in W^{2,2}(I,M)$ and therefore it is $C^1$ on $I$. It implies that $\dot{\gamma}_{+}(0)=\dot{\gamma}_{-}(1)$.
		
		To prove the converse statement, let $\gamma$ be a closed weak geodesic in $S$ and the sectional curvatures of $M$ are nonpositive in a neighborhood of $\im(\gamma)$. Then using Theorem \ref{weaklocal}, $\gamma$ is locally minimizing. We now suppose that $W$ is a tubular  neighborhood of $\im(\gamma)$ with compact closure and with radius smaller than $\bar{r}$, where $\bar{r}=\min_{t\in I}r\left(\gamma(t)\right)$ and  $\se\leq 0$ on $W$.  Let $\eta\in \mathcal{C}$ be such that $\im(\eta)\subset W$, hence  $c\left(\gamma,\eta\right)\equiv 2$ on $I$, where $c\left(\gamma,\eta\right)(t)=2\sqrt{\Delta}\:d\left(\gamma(t),\eta(t)\right)
		\cot\left(\sqrt{\Delta}\:d\left(\gamma(t),\eta(t)\right)\right)$ for all $t\in I$. Therefore similar to the proof of Theorem \ref{weaklocal} we have
		\begin{equation}\label{ineqq1}
			\frac{1}{2}\norm \dot{\eta}\norm^2-\frac{1}{2}\norm \dot{\gamma}\norm^2\geq \frac{1}{2}\norm L_{\eta,\gamma}\dot{\eta}-\dot{\gamma}\norm^2+\frac{d}{dt}\left\li \dot{\gamma}, \exp_{\gamma}^{-1}\eta\right\ri
		\end{equation}
		\[
		-  \left\li D_t\dot{\gamma}, \exp_{\gamma}^{-1}\eta\right\ri-\frac{1}{2}|R|_{\infty} d^2\left(\gamma, \eta\right)\norm \dot{\gamma}\norm \norm \dot{\eta}\norm,
		\]
		where $R$ denotes the curvature tensor on $M$.
		Hence by integrating from both side of \eqref{ineqq1} and noting that $D_t\dot{\gamma}\in N^P_S\left(\gamma\right),\ a.\ e.$ and $\dot{\gamma}_{+}(0)=\dot{\gamma}_{-}(1)$, we have
		\begin{equation}\label{ineqq2}
			\frac{1}{2}\int_{0}^{1}\norm \dot{\eta}(t)\norm^2dt -  \frac{1}{2}\int_{0}^{1}\norm \dot{\gamma}(t)\norm^2dt\geq
		\end{equation}
		\[
		\frac{1}{2}\int_{0}^{1}\norm L_{\eta,\gamma}\dot{\eta}-\dot{\gamma}\norm^2dt-\left(\bar{\varphi}\int_{0}^{1}\norm D_t\dot{\gamma}\norm dt+\frac{1}{2}|R|_{\infty} \norm \dot{\gamma}\norm_{L^2} \norm \dot{\eta}\norm_{L^2} \right)d^2_{\infty}\left(\gamma,\eta\right).
		\]

		On the other hand using Lemma \ref{dinfty} we have
		\begin{equation}\label{ineq4}
			d_{\infty}^2\left(\gamma,\eta\right)\leq h\int_{0}^{1}\norm L_{\eta,\gamma}\dot{\eta}-\dot{\gamma}\norm^2dt,
		\end{equation}
		for two closed curves $\gamma,\eta\in H^1(I,M)$ with the properties that $\eta=\gamma$ on $[h,1-h]$ for some $h>0$ and  $d_{\infty}\left(\gamma, \eta\right)< \bar{r}$. Indeed, it suffices to define the new curves $\tilde{\gamma},\tilde{\eta}\in H^1(I,M)$ by
		\[
		\tilde{\gamma}(t):=\left\{
		\begin{array}{ll}
			\gamma\left(t+\frac{1}{2}\right) & 0\leq t\leq \frac{1}{2}  \\
			\gamma\left(t-\frac{1}{2}\right) & \frac{1}{2}\leq t\leq 1,
		\end{array}
		\right.
		\]
		and
		\[
		\tilde{\eta}(t):=\left\{
		\begin{array}{ll}
			\eta\left(t+\frac{1}{2}\right) & 0\leq t\leq \frac{1}{2}  \\
			\eta\left(t-\frac{1}{2}\right) & \frac{1}{2}\leq t\leq 1.
		\end{array}
		\right.
		\]
		Hence $\tilde{\eta}=\tilde{\gamma}$ on $I\setminus \left(\frac{1}{2}-h,\frac{1}{2}+h\right)$. Then using Lemma \ref{dinfty} and changing the variables, we get the inequality \eqref{ineq4}.
		
		We now put
		\[
		\theta:=\bar{\varphi}\int_{0}^{1}\norm D_t\dot{\gamma}(t)\norm dt+\frac{1}{2}C|R|_{\infty}\norm \dot{\gamma}\norm_{L^2},
		\]
		where $C>0$ is such that $C>\norm \dot{\gamma}\norm_{L^2}$ and choose $h>0$ small enough such that $h<\frac{1}{2\theta}$. Let $\eta\in \mathcal{C}\setminus\{\gamma\}$ be such that $\eta=\gamma$ on $ \left[h,1-h\right]$, $\norm \dot{\eta}\norm_{L^2}<C$ and $\im(\eta)\subset W$. Then using \eqref{ineq4} we obtain that
		\[
		\frac{1}{2}\int_{0}^{1}\norm \dot{\eta}(t)\norm^2dt- \frac{1}{2}\int_{0}^{1}\norm \dot{\gamma}(t)\norm^2dt\geq
		\left(\frac{1}{2}-h\theta\right)\int_{0}^{1}\norm L_{\eta,\gamma}\dot{\eta}-\dot{\gamma}\norm^2dt.
		\]
		Moreover, by \eqref{ineq4} we have $\int_{0}^{1}\norm L_{\eta,\gamma}\dot{\eta}-\dot{\gamma}\norm^2dt=0$ if and only if $\eta=\gamma$. It follows that
		\begin{equation*}
			\frac{1}{2}\int_{0}^{1}\norm \dot{\eta}(t)\norm^2dt > \frac{1}{2}\int_{0}^{1}\norm \dot{\gamma}(t)\norm^2dt.
		\end{equation*}
		
		We note that if $\im(\eta)\nsubseteq W$, then $\mathcal{L}(\gamma)<\mathcal{L}(\eta)$ and it completes the proof.
		
	\end{proof}



\section*{Statements and Declarations}

The first-named author and the second-named author declare they have no financial interests. The third-named author was supported by the Iran National Science Foundation (INSF) under project No.4002602.
	

	
\end{document}